\documentclass[11pt]{article}
\usepackage{amsfonts}
\usepackage{bbm}
\usepackage{mathrsfs}
\usepackage{amssymb}

\def\hang{\hangindent\parindent}
\def\textindent#1{\indent\llap{#1\enspace}\ignorespaces}
\def\re{\par\hang\textindent}

\title{Note on (De)homogenized Gr\"obner Bases\thanks{Project supported by the National Natural
Science Foundation of China (10571038).}}

\author{Huishi Li\thanks{e-mail: huishipp@yahoo.com}\\
{\small Department of Applied Mathematics}\\
{\small College of Information Science and Technology}\\
{\small Hainan University}\\
{\small  Haikou 570228, China}}

\date{}

\begin{document}
\maketitle
\begin{center}
\begin{minipage}{120mm}
{\small {\bf Abstract.} By employing the (de)homogenization
technique in a relatively extensive setting, this note studies in
detail the relation between non-homogeneous Gr\"obner bases and
homogeneous Gr\"obner bases. As a consequence,  a general principle
of computing Gr\"obner bases (for an ideal and its homogenization
ideal) by passing to homogenized generators is clarified
systematically. The obtained results improve and strengthen the work
of [LWZ], [Li1], [Li2], [Li3], and very recent [SL] concerning the
same topic. }
\end{minipage}\end{center}{\parindent=0pt\vskip 6pt
{\bf 2000 Mathematics Classification} Primary 16W70; Secondary 68W30
(16Z05).\vskip 6pt {\bf Key words} Graded algebra, graded ideal,
Gr\"obner basis, homogenization, dehomogenization }

\def\QED{\hfill{$\Box$}} \def\NZ{\mathbb{N}}
\def \r{\rightarrow}

\def\normalbaselines{\baselineskip 24pt\lineskip 4pt\lineskiplimit 4pt}
\def\mapdown#1{\llap{$\vcenter {\hbox {$\scriptstyle #1$}}$}
                                \Bigg\downarrow}
\def\mapdownr#1{\Bigg\downarrow\rlap{$\vcenter{\hbox
                                    {$\scriptstyle #1$}}$}}
\def\mapright#1#2{\smash{\mathop{\longrightarrow}\limits^{#1}_{#2}}}
\def\mapleft#1#2{\smash{\mathop{\longleftarrow}\limits^{#1}_{#2}}}
\def\mapup#1{\Bigg\uparrow\rlap{$\vcenter {\hbox  {$\scriptstyle #1$}}$}}
\def\mapupl#1{\llap{$\vcenter {\hbox {$\scriptstyle #1$}}$}
                                      \Bigg\uparrow}
\def\v5{\vskip .5truecm}
\def\T#1{\widetilde #1}
\def\OV#1{\overline {#1}}
\def\hang{\hangindent\parindent}
\def\textindent#1{\indent\llap{#1\enspace}\ignorespaces}
\def\item{\par\hang\textindent}

\def\LH{{\bf LH}}\def\LM{{\bf LM}}\def\LT{{\bf
LT}}\def\KX{K\langle X\rangle} \def\KS{K\langle X\rangle}
\def\B{{\cal B}} \def\LC{{\bf LC}} \def\G{{\cal G}} \def\FRAC#1#2{\displaystyle{\frac{#1}{#2}}}
\def\SUM^#1_#2{\displaystyle{\sum^{#1}_{#2}}} \def\O{{\cal O}}  \def\J{{\bf J}}
{\parindent=0pt\vskip 1truecm

In the computational Gr\"obner basis theory, though it is a
well-known fact that by virtue of both the structural advantage
(mainly the degree-truncated structure) and the computational
advantage (mainly the degree-preserving fast ordering), most of the
popularly used commutative and noncommutative Gr\"obner basis
algorithms produce Gr\"obner bases by homogenizing generators first,
it seems that in both the commutative and noncommutative case a
general principle of computing Gr\"obner bases (for an ideal and its
homogenization ideal) by passing to homogenized generators is still
missing.}\par
Let $K$ be a field, and let $\NZ$ be the additive monoid of
nonnegative integers. Recall from [Li2] that if
$R=\oplus_{p\in\NZ}R_p$ is an $\NZ$-graded $K$-algebra with an
admissible system $(\B ,\prec)$, in which $\B$ is a {\it skew
multiplicative $K$-basis} of $R$ consisting of $\NZ$-{\it
homogeneous elements} (i.e., $u,v\in\B$ implies that $u$, $v$ are
homogeneous elements, $uv=0$ or $uv=\lambda w$ for some nonzero
$\lambda\in K$ and $w\in\B$), and $\prec$ is a monomial ordering on
$\B$, then, theoretically every (two-sided) ideal $I$ of $R$ has a
Gr\"obner basis $\G$ in the sense that if $f\in I$ and $f\ne 0$ then
there is some $g\in\G$ such that $\LM (g)|\LM (f)$, where $\LM (~)$
denotes taking the leading monomial of elements in $R$, in
particular, every graded ideal of $R$ has a {\it homogeneous
Gr\"obner basis}, i.e., a Gr\"obner basis consisting of
$\NZ$-homogeneous elements. Typical examples of such algebras
include commutative polynomial $K$-algebra, noncommutative free
$K$-algebra, path algebra over $K$, the coordinate algebra of a
quantum affine $n$-space over $K$, and exterior $K$-algebra (cf.
[Bu], [BW], [Gr], [HT], [KRW], [Mor]).
 In this note, we employ the
(de)homogenization technique to study in detail the relation between
Gr\"obner bases in $R$ and homogeneous Gr\"obner bases in the
polynomial ring $R[t]$, respectively the relation between Gr\"obner
bases in the free algebra $\KS =K\langle X_1,...,X_n\rangle$ and
homogeneous Gr\"obner bases in the free algebra $K\langle X,T\rangle
=K\langle X_1,...,X_n,T\rangle$.  As a consequence, this makes a
solid theoretical foundation for us to demonstrate, by passing to
the graded ideal $\langle S^*\rangle$, how to obtain a Gr\"obner
basis for the ideal $I=\langle S\rangle$ generated by a subset
$S\subset R$ and hence a homogeneous Gr\"obner basis for the central
homogenization ideal $\langle I^*\rangle$ of $I$ in $R[t]$ with
respect to $t$;  respectively, this enables us to demonstrate, by
passing to the graded ideal $\langle~\T S~\rangle$, how to obtain a
Gr\"obner basis for the ideal $I=\langle S\rangle$ generated by a
subset $S\subset \KS$ and hence a homogeneous Gr\"obner basis for
the non-central homogenization ideal $\langle~\T I~\rangle$ of $I$
in $K\langle X,T\rangle$ with respect to $T$. The obtained results
improve and strengthen the work of [LWZ], [Li1], [Li3], [Li3], and
very recent [SL] concerning the same topic.  \v5

Algebras considered in this paper are associative algebras with
multiplicative identity 1. Unless otherwise stated, ideals
considered are meant two-sided ideals. If $S$ is a nonempty subset
of an algebra, then we use $\langle S\rangle$ to denote the
two-sided ideal generated by $S$. Moreover, if $K$ is a field, then
we write $K^*=K-\{ 0\}$.\v5

\section*{1. Central (De)homogenized Gr\"obner Bases}
Let $R=\oplus_{p\in\NZ}R_p$ be an arbitrary $\NZ$-graded
$K$-algebra, and let $R[t]$ be the polynomial ring in the commuting
variable $t$ over $R$. Then $R[t]$ has the  mixed $\NZ$-gradation,
that is, $R[t]=\oplus_{p\in\NZ}R[t]_p$ is an $\NZ$-graded algebra
with the degree-$p$ homogeneous part
$$R[t]_p=\left\{\left. \displaystyle{\sum_{i+j=p}}F_it^j~\right |~F_i\in R_i,
~j\ge 0
 \right \} ,\quad p\in \mathbb{N} .$$
Considering the onto ring homomorphism $\phi$: $R[t]\rightarrow R $
defined by $\phi (t)=1$, then for each $f\in R$, there exists a
homogeneous element $F\in R[t]_p$, for some $p$, such that $\phi
(F)=f$. More precisely, if $f=f_p+f_{p-1}+ \cdots +f_{p-s}$ with
$f_p\in R_p$,  $f_{p-j}\in R_{p-j}$ and $f_p\ne 0$, then
$$f^*=f_p+tf_{p-1}+\cdots +t^sf_{p-s}$$ is a homogeneous element of
degree $p$ in $R[t]_p$ satisfying $\phi (f^*)=f$. We call the
homogeneous element $f^*$ obtained this way the {\it central
homogenization} of $f$ with respect to $t$ (for the reason that $t$
is in the center of $R[t]$). On the other hand, for an element $F\in
R[t]$, we write
$$F_*=\phi (F)$$ and call it the {\it central dehomogenization} of
$F$ with respect to $t$ (again for the reason that $t$ is in the
center of $R[t]$). Hence, if $I$ is an ideal $R$, then we write
$I^*=\{ f^*~|~f\in I\}$ and call the $\NZ$-graded ideal $\langle
I^*\rangle$ generated by $I^*$ the {\it central homogenization
ideal} of $I$ in $R[t]$ with respect to $t$; and if $J$ is an ideal
of $R[t]$, then since $\phi$ is a ring epimorphism, $\phi (J)$ is an
ideal of $R$, so we write $J_*$ for $\phi (J)=\{H_*=\phi (H)~|~H\in
J\}$ and call it the {\it central dehomogenization ideal} of $J$ in
$R$ with respect to $t$. Consequently, henceforth we will also use
the notation
$({\bf\emph{J}}_*)^*=\{(h_*)^*\mid{h\in{\bf\emph{J}}}\}$. \par
Since each $f\in R$ has a unique decomposition by $\NZ$-homogeneous
elements, if if $f=f_p+f_{p-1}+ \cdots +f_{p-s}$ with $f_p\in R_p$,
$f_{p-j}\in R_{p-j}$ and $f_p\ne 0$, then we call $f_p$ the
$\NZ$-leading homogeneous element of $f$, denoted $\LH (f)=f_p$.
Similarly, for an element $F\in R[t]$, the $\NZ$-leading homogeneous
element $\LH (F)$ is defined with respect to the mixed
$\NZ$-gradation of $R[t]$. {\parindent=0pt\v5

{\bf 1.1. Lemma} With every definition and notation made above, the
following statements hold.\par (i) For $F,G\in R[t]$,
$(F+G)_*=F_*+G_*$, $(FG)_*=F_*G_*$.\par (ii) For any $f\in R$,
$(f^*)_*=f$.\par (iii) If $F\in R[t]_p$ and if $(F_*)^*\in R[t]_q$,
then $p\ge q$ and $t^r(F_*)^*=F$ with $r=p-q$.\par (iv) If $f,g\in
R$ are such that $fg$ has nonzero $\LH (fg)\in R_m$ and $f^*g^*$ has
nonzero $\LH (f^*g^*)\in R[t]_q$, then $f^*g^*=t^k(fg)^*$ with
$k=q-m$. \par (v) If $f,g\in R$ with nonzero $\LH (f)\in R_p$ and
nonzero $\LH (g)\in R_q$, then $(f+g)^*=f^*+g^*$ in case $p=q$, and
$(f+g)^*=f^*+t^{\ell}g^*$ in case $p>q$, where $\ell =p-q$.\par (vi)
If $I$ is a two-sided ideal of $R$, then each homogeneous element
$F\in \langle I^*\rangle$ is of the form $t^rf^*$ for some $r\in\NZ$
and $f\in I$.\par (vii) If $J$ is a graded ideal of $R[t]$, then for
each $h\in J_*$ there is some homogeneous element $F\in J$ such that
$F_*=h$.\vskip 6pt {\bf Proof} Exercise. \QED}\v5

Suppose that the $\NZ$-graded $K$-algebra $R=\oplus_{p\in\NZ}R_p$
has an admissible system $(\B ,\prec_{gr})$, where $\B$ is a skew
multiplicative $K$-basis of $R$ consisting of $\NZ$-{\it homogeneous
elements}, and $\prec_{gr}$ is an $\NZ$-graded monomial ordering on
$\B$, i.e., $R$ has a Gr\"obner basis theory. Consider the mixed
$\NZ$-gradation of $R[t]$ and the $K$-basis $\B^*=\{ t^rw~|~w\in \B
,~r\in\NZ\}$ of $R[t]$. Since $\B^*$ is obviously a skew
multiplicative $K$-basis for $R[t]$, the $\NZ$-graded monomial
ordering $\prec_{gr}$ on $\B$ extends to a monomial ordering on
$\B^*$, denoted $\prec_{t\hbox{-}gr}$, as follows:
$$t^{r_1}w_1\prec_{t\hbox{-}gr}t^{r_2}w_2~\hbox{if~and ~only~if~} w_1\prec_{gr}
w_2,~\hbox{or}~w_1=w_2~\hbox{and}~r_1<r_2.$$ Thus $R[t]$ holds a
Gr\"obner basis theory  with respect to the admissible system
$(\B^*,\prec_{t\hbox{-}gr})$.\par
As usual we call elements in $\B$ and $\B^*$ monomials and use $\LM
(~)$ to denote taking the leading monomial of elements with respect
to the given monomial ordering.\v5

It follows from the definition of $\prec_{t\hbox{-}gr}$ that
$t^r\prec_{t\hbox{-}gr}w$ for all integers $r>0$ and all $w\in\B -\{
1\}$ (if $\B$ contains the identity element 1 of $R$). Hence
$\prec_{t\hbox{-}gr}$ {\it is not a graded monomial ordering} on
$\B^*$. But  noticing that $\prec_{gr}$ is an $\NZ$-graded monomial
ordering on $\B$, under taking the $N$-leading homogeneous element
and central (de)homogenization the leading monomials behave
harmonically as described in the next lemma. {\parindent=0pt\v5

{\bf 1.2. Lemma} With notation given above, the following statements
hold.\par (i) If $f\in R$, then $$\LM (f)=\LM (\LH
(f))~\hbox{w.r.t.}~\prec_{gr}~\hbox{on}~\B .$$\par (ii) If $f\in R$,
then
$$\LM (f^*)=\LM (f)~\hbox{w.r.t.}~\prec_{t\hbox{-}gr}~\hbox{on}~\B^*.$$\par
(iii) If $F$ is a nonzero homogeneous element of $R[t]$, then
$$\LM (F_*)=\LM (F)_*~\hbox{w.r.t.}~\prec_{gr}~\hbox{on}~\B .$$\par
{\bf Proof} The proof of (i) and (ii) is an easy exercise. To prove
(iii), let $F\in R[t]_p$ be a nonzero homogeneous element of degree
$p$, say $$F=\lambda t^rw+\lambda_1t^{r_1}w_1+\cdots
+\lambda_st^{r_s}w_s,$$ where $\lambda,\lambda_i,\in K^*$,
$r,r_i\in\NZ$,  $w,w_i\in\B$, such that $\LM (F)=t^rw$. Since $\B$
consists of $\NZ$-homogeneous elements and $R[t]$ has the mixed
$\NZ$-gradation by the previously fixed assumption, we have
$d(t^rw)=d(t^{r_i}w_i)=p$, $1\le i\le s$. Thus $w=w_i$ will imply
$r=r_i$ and thereby $t^rw=t^{r_i}w_i$. So we may assume that $w\ne
w_i$, $1\le i\le s$. Then it follows from the definition of
$\prec_{t\hbox{-}gr}$ that $w_i\prec_{gr}w$ and $r\le r_i$, $1\le
i\le s$. Therefore $\LM (F_*)=w=\LM (F)_*$, as desired.\QED}\v5

The next result is a generalization of ([LWZ] Theorem
2.3.2).{\parindent=0pt\v5

{\bf 1.3. Theorem} With notions and notations as fixed before, let
$I=\langle\G\rangle$ be the ideal of $R$ generated by a subset $\G$,
and $\langle I^*\rangle$ the central homogeneization ideal of $I$ in
$R[t]$ with respect to $t$. The following two statements are
equivalent.\par (i) $\G$ is a Gr\"obner basis for $I$ in $R$ with
respect to the admissible system $(\B ,\prec_{gr})$;\par (ii)
$\G^*=\{ g^*~|~g\in\G\}$  is a Gr\"obner basis for $\langle
I^*\rangle$ in $R[t]$ with respect to the admissible system
$(\B^*,\prec_{t\hbox{-}gr})$. \vskip 6pt {\bf Proof} In proving the
equivalence below, without specific indication we shall use (i) and
(ii) of Lemma 1.2 wherever it is needed.
\par (i) $\Rightarrow$ (ii) First note that $\LM (\G^*)\subset\LM
(I^*)$. We have to prove that $\LM (\G^*)$ generates $\langle\LM
(I^*)\rangle$ in order to see that $\G^*$ is a Gr\"obner basis for
$\langle I^*\rangle$. If $F\in \langle I^*\rangle$, then since $\LM
(F)=\LM (\LH (F))$, we may assume, without loss of generality, that
$F$ is a homogeneous element. So, by Lemma 1.1(vi) we have
$F=t^rf^*$ for some $f\in I$. It follows from the equality $\LM
(f^*)=\LM (f)$ that
$$\LM (F)=t^r\LM (f^*)=t^r\LM (f).$$
Since $\G$ is a Gr\"obner basis for $I$, $\LM (f)=\lambda v\LM
(g_i)w$ for some $\lambda\in K^*$, $g_i\in\G$, and $v,w\in\B$. Thus,
$$\LM (F)=t^r\LM (f)=\lambda t^rv\LM (g_i^*)w\in\langle\LM (\G^*)\rangle .$$
This shows that $\langle\LM (\langle I^*\rangle )\rangle =\langle\LM
(\G^*)\rangle$, as desired.\par (ii) $\Rightarrow$ (i) Suppose
$\G^*$ is a Gr\"obner basis for the homogenization ideal $\langle
I^*\rangle$ of $I$ in $R[t]$.  Let $f\in I$. Then $\LM (f^*)=\lambda
v\LM (g_i^*)w$ for some $\lambda\in K^*$, $v,w\in\B^*$ and
$g_i^*\in\G^*$. Since $\LM (f)=\LM (f^*)$, it follows that
$$\LM (f)=\lambda v_*\LM (g_i)w_*\in \langle\LM (\G )\rangle .$$
This shows that $\langle\LM (I)\rangle =\langle\LM (\G )\rangle$,
i.e., $\G$ is a Gr\"obner basis for $I$ in $R$.\QED}\v5

We call the Gr\"obner basis $\G^*$ obtained in Theorem 1.3 the {\it
central homogenization of} $\G$ in $R[t]$ with respect to $t$, or
$\G^*$ is a {\it central homogenized Gr\"obner basis} with respect
to $t$.\v5

By Lemma 1.2  and Theorem 1.3, we have immediately the following
corollary. {\parindent=0pt\v5

{\bf 1.4. Corollary} Let $I$ be an arbitrary ideal of $R$.  With
notation as before, if $\G$ is a Gr\"obner basis of $I$ with respect
to the data $(\B ,\prec_{gr})$, then, with respect to the data
$(\B^*, \prec_{t\hbox{-}gr})$ we have
$$\B^*-\langle\LM (\G^*)\rangle =\{ t^rw~|~w\in\B -\langle\LM (\G )\rangle,~r\in\NZ\} ,$$
that is, the set $N(\langle I^*\rangle )$ of normal monomials (mod
$\langle I^*\rangle )$ in $\B^*$ is determined by the set $N(I)$ of
normal monomials (mod $I$) in $\B$. Hence, the  algebra
$R[t]/\langle I^*\rangle =R[t]/\langle\G^*\rangle$ has the $K$-basis
$$\OV{N(\langle I^*\rangle )}=\{ \OV{t^rw}~|~w\in
N(I),~r\in\NZ\}.$$}\QED\v5

Theoretically we may also obtain a Gr\"obner basis for an ideal $I$
of $R$ by dehomogenizing a homogeneous Gr\"obner basis of the ideal
$\langle I^*\rangle\subset R[t]$. Below we give a more general
approach to this assertion.{\parindent=0pt \v5

{\bf 1.5. Theorem} Let $J$ be a graded ideal of $R[t]$. If
$\mathscr{G}$ is a homogeneous Gr\"obner basis of $J$ with respect
to the data $(\B^*, \prec_{t\hbox{-}gr})$, then $\mathscr{G}_*=\{
G_*~|~G\in\mathscr{G}\}$ is a Gr\"obner basis for the ideal $J_*$ in
$R$ with respect to the data $(\B ,\prec_{gr})$.\vskip 6pt {\bf
Proof} If $\mathscr{G}$ is a Gr\"obner basis of $J$, then
$\mathscr{G}$ generates $J$ and hence $\mathscr{G}_*=\phi
(\mathscr{G} )$ generates $J_*=\phi (J)$. For a nonzero $f\in J_*$,
by Lemma 1.1(vii), there exists a homogeneous element $H\in J$ such
that $H_*=f$. It follows from Lemma 1.2 that
$$\LM (f)=\LM (f^*)=\LM ((H_*)^*).\leqno{(1)}$$
On the other hand, there exists some $G\in\mathscr{G}$ such that
$\LM (G)|\LM (H)$, i.e., $$\LM (H)=\lambda t^{r_1}w\LM
(G)t^{r_2}v\leqno{(2)}$$ for some $\lambda\in K^*$, $r_1,r_2\in\NZ$,
$w,v\in\B$. But by Lemma 1.1(iii) we also have $t^r(H_*)^*=H$ for
some $r\in\NZ$, and hence
$$\LM (H)=\LM (t^r(H_*)^*)=t^r\LM (H_*)^*).\leqno{(3)}$$
So, $(1)+(2)+(3)$ yields
$$\begin{array}{rcl} \lambda t^{r_1+r_2}w\LM (G)v&=&\LM (H)\\
&=&t^r\LM ((H_*)^*)\\
&=&t^r\LM (f).\end{array}$$ Taking the central dehomogenization for
the above equality, by Lemma 1.2(iii) we obtain $$\lambda w\LM
(G_*)v=\lambda w\LM (G)_*v=\LM (f).$$ This shows that $\LM (G_*)|\LM
(f)$. Therefore, $\mathscr{G}_*$ is a Gr\"obner basis for
$J_*$.\QED}\v5

We call the Gr\"obner basis $\mathscr{G}_*$ obtained in Theorem 1.5
the {\it central dehomogenization of} $\mathscr{G}$ in $R$ with
respect to $t$, or $\mathscr{G}_*$ is a {\it central dehomogenized
Gr\"obner basis} with respect to $t$.{\parindent=0pt\v5

{\bf 1.6. Corollary} Let $I$ be an ideal of $R$. If $\mathscr{G}$ is
a homogeneous Gr\"obner basis of $\langle I^*\rangle$ in $R[t]$ with
respect to the data $(\B^*, \prec_{t\hbox{-}gr})$, then
$\mathscr{G}_*=\{ g_*~|~g\in\mathscr{G}\}$ is a Gr\"obner basis for
$I$ in $R$ with respect to the data $(B ,\prec_{gr})$. Moreover,  if
$I$ is generated by the subset $F$ and $F^*\subset\mathscr{G}$, then
$F\subset\mathscr{G}_*$. \vskip 6pt {\bf Proof} Put $J=\langle
I^*\rangle$. Then since $J_*=I$, it follows from Theorem 1.5 that if
$\mathscr{G}$ is a homogeneous Gr\"obner basis of $J$ then
$\mathscr{G}_*$ is a Gr\"obner basis for $I$. The second assertion
of the theorem is clear by Lemma 1.1(ii).\QED}\v5

Let $S$ be a nonempty subset of $R$ and $I=\langle S\rangle$ the
ideal generated by $S$. Then, with $S^*=\{ f^*~|~f\in S\}$, in
general $\langle S^*\rangle\subsetneq\langle I^*\rangle$ in $R[t]$
(for instance, consider $S=\{ y^3-x-y,~y^2+1\}$ in the commutative
polynomial ring $K[x,y]$ and the homogenization in $K[x,y,t]$ with
respect to $t$). So, from both a practical and a computational
viewpoint, it is the right place to set up the procedure of getting
a Gr\"obner basis for $I$ and hence a Gr\"obner basis for $\langle
I^*\rangle$ by producing a homogeneous Gr\"obner basis of the graded
ideal $\langle S^*\rangle$.{\parindent=0pt\v5

{\bf 1.7. Proposition} Let $I=\langle S\rangle$ be the ideal of $R$
as fixed above. Suppose that Gr\"obner bases are algorithmically
computable in $R$ and hence in $R[t]$. Then a Gr\"obner basis for
$I$ and a homogeneous Gr\"obner basis for $\langle I^*\rangle$ may
be obtained by implementing the following procedure:\par {\bf Step
1.}  Starting with the initial subset $S^*=\{ f^*~|~f\in S\}$,
compute a homogeneous Gr\"obner basis $\mathscr{G}$ for the graded
ideal $\langle S^*\rangle$ of $R[t]$.\par {\bf Step 2.} Noticing
$\langle S^*\rangle_*=I$, use Theorem 1.5 and dehomogenize
$\mathscr{G}$ with respect to $t$ in order to obtain the Gr\"obner
basis $\mathscr{G}_*$ for $I$.\par {\bf Step 3.} Use Theorem 1.3 and
homogenize $\mathscr{G}_*$ with respect to $t$ in order to obtain
the homogeneous Gr\"obner basis $(\mathscr{G}_*)^*$ for the graded
ideal $\langle I^*\rangle$.\par\QED}\v5
\def\KXT{K\langle X, T\rangle}

\section*{2. Non-central (De)homogenized Gr\"obner Bases}
In a similar way, we proceed now to consider the free algebra
$\KS=K\langle X_1,...,X_n\rangle$ of $n$ generators as well as the
free algebra $K\langle X, T\rangle =K\langle X_1,...,X_n,T\rangle$
of $n+1$ generators, and demonstrate how Gr\"obner bases in $\KS$
are related with homogeneous Gr\"obner bases in $\KXT$ if the
non-central (de)homogenization with respect to $T$ is employed. \v5

Let $\KS$ be equipped with a fixed {\it weight $\NZ$-gradation}, say
each $X_i$ has degree $n_i>0$, $1\le i\le n$. Assigning to $T$ the
degree 1 in $\KXT$ and using the same weight $n_i$ for each  $X_i$
as in $\KS$, we get the weight $\NZ$-gradation of $K\langle
X,T\rangle$ which extends the weight $\NZ$-gradation of $\KS$. Let
$\B$ and $\T{\B}$ denote the standard $K$-bases of $\KS$ and $\KXT$
respectively. To be convenient we use lowercase letters $w,u,v,...$
to denote monomials in $\B$ as before, but use capitals $W,U,V,...$
to denote monomials in $\T B$.\par In what follows, we fix an
admissible system $(\B, \prec_{gr})$ for $\KS$, where $\prec_{gr}$
is an $\NZ$-{\it graded lexicographic ordering} on $\B$ with respect
to the fixed weight $\NZ$-gradation of $\KS$, such that
$$X_{i_1}\prec_{gr}X_{i_2}\prec_{gr}\cdots\prec_{gr}X_{i_n}.$$ Then it is not
difficult to see that $\prec_{gr}$ can be extended to an $\NZ$-{\it
graded lexicographic ordering} $\prec_{_{T\hbox{-}gr}}$ on $\T{\B}$
with respect to the fixed weight $\NZ$-gradation  of $\KXT$, such
that
$$T\prec_{_{T\hbox{-}gr}} X_{i_1}\prec_{_{T\hbox{-}gr}}X_{i_2}\prec_{_{T\hbox{-}gr}}\cdots
\prec_{_{T\hbox{-}gr}}X_{i_n},$$ and thus we get the admissible
system $(\T{\B}, \prec_{_{T\hbox{-}gr}})$ for $\KXT$. With respect
to $\prec_{gr}$ and $ \prec_{_{T\hbox{-}gr}}$ we use $\LM (~)$ to
denote taking the leading monomial of elements in $\KS$ and $\KXT$
respectively. \v5

Consider the fixed $\NZ$-graded structures $\KS
=\oplus_{p\in\NZ}\KS_p$, $\KXT =\oplus_{p\in\NZ}\KXT_p$, and the
ring epimorphism
$$\psi :~\KXT~\longrightarrow~\KS$$ defined by
$\psi (X_i)=X_i$ and $\psi (T)=1$. Then each $f\in \KS$ is the image
of some homogeneous element in $\KXT$. More precisely, if
$f=f_{p}+f_{p-1}+\cdots+f_{p-s}$ with $f_{p}\in \KS_p$, $f_{p-j}\in
\KS_{p-j}$ and $f_{p}\neq0$, then $$\T
f=f_{p}+Tf_{p-1}+\cdots+T^sf_{p-s}$$ is a homogeneous element of
degree $p$ in $\KXT_p$ such that $\psi (\T f)=f$. We call the
homogeneous element $\T f$ obtained this way the {\it non-central
homogenization} of $f$ with respect to $T$ (for the reason that $T$
is not a commuting variable). On the other hand, for $F\in \KXT$, we
write $$F_{\sim}=\psi(F)$$ and call $F_{\sim}$ the {\it non-central
dehomogenization} of $F$ with respect to $T$ (again for the reason
that $T$ is not a commuting variable). Furthermore, if $I=\langle
S\rangle$ is the ideal of $\KS$ generated by a subset $S$, then we
write
$$ \begin{array}{l}
\T S=\{ \T f~|~f\in S\}\cup\{X_iT-TX_i~|~1\le i\le n\},\\  \T I=\{
\T f~|~f\in I\}\cup\{X_iT-TX_i~|~1\le i\le n\},\end{array}$$ and
call the graded ideal $\langle~\T I~\rangle$ generated by $\T I$ the
{\it non-central homogenization ideal} of $I$ in $\KXT$ with respect
to $T$; while if $J$ is an ideal of $\KXT$, then since $\psi$ is a
surjective ring homomorphism, $\psi (J)$ is an ideal of $\KS$, so we
write $J_{\sim}$ for $\psi (J)=\{H_{\sim}|H\in J\}$ and call it the
{\it non-central dehomogenization ideal} of $J$ in $\KS$ with
respect to $T$. Consequently, henceforth we will also use the
notation
$$(J_{\sim})^{\sim}=\{(h_{\sim})^{\sim}~|~h\in J\}\cup\{
X_iT-TX_i~|~1\le i\le n\} .$$\par It is straightforward to check
that with resspect to the data $(\T{\B},\prec_{_{T\hbox{-}gr}})$,
the subset $\{ X_iT-TX_i~|~1\le i\le n\}$ of $\KXT$ forms a
homogeneous Gr\"obner basis with $\LM (X_iT-TX_i)=X_iT$, $1\le i\le
n$. In the latter discussion we will freely use this fact without
extra indication.{\parindent=0pt\v5

{\bf 2.1. Lemma} With notation as fixed before, the following
properties hold.\par (i) If $F,G\in \KXT$, then
$(F+G)_{\sim}=F_{\sim}+G_{\sim}$,
$(FG)_{\sim}=F_{\sim}G_{\sim}$.\par (ii) For each nonzero $f\in
\KS$, $(\T f)_{\sim}=f$.\par (iii) Let $\mathscr{C}$ be the graded
ideal of $\KXT$ generated by $\{ X_iT-TX_i~|~1\le i\le n\}$. If
$F\in \KXT_p$, then there exists an $L\in \mathscr{C}$ and a unique
homogeneous element of the form $H=\sum_i\lambda_iT^{r_i}w_i$, where
$\lambda_i\in K^*$, $w_i\in\B$, such that $F=L+H$; moreover there is
some $r\in\NZ$ such that $T^r(H_{\sim})^{\sim}=H$, and hence
$F=L+T^r(F_{\sim})^{\sim}$.
\par
(iv) Let $\mathscr{C}$ be as in (iii) above.  If $I$ is an ideal of
$\KS$, $F\in\langle ~\T I~\rangle$ is a homogeneous element, then
there exist some $L\in \mathscr{C}$, $f\in I$ and $r\in\NZ$ such
that $F=L+T^r\T f$.\par (v) If $J$ is a graded ideal of $\KXT$ and
$\{ X_iT-TX_i~|~1\le i\le n\}\subset J$, then for each nonzero $h\in
J_{\sim}$, there exists a homogeneous element
$H=\sum_i\lambda_iT^{r_i}w_i\in J$, where $\lambda_i\in K^*$,
$r_i\in\NZ$, and $w_i\in\B$, such that for some $r\in\NZ$,
$T^r(H_{\sim})^{\sim}=H$ and $H_{\sim}=h$. \vskip 6pt {\bf Proof}
(i) and (ii) follow from the definitions of non-central
homogenization and  non-central dehomogenization directly.\par (iii)
Since the subset $\{ X_iT-TX_i~|~1\le i\le n\}$ is a Gr\"obner basis
in  $\KXT$ with respect to $(\T{\B},\prec_{_{T\hbox{-}gr}})$, such
that  $\LM (X_iT-TX_i)=X_iT$, $1\le i\le n$,  if $F\in\KXT_p$, then
the division of $F$ by this subset yields $F=L+H$, where $L\in
\mathscr{C}$, and $H=\sum_i\lambda_iT^{r_i}w_i$ is the unique
remainder with $\lambda_i\in K^*$, $w_i\in\B$, in which each
monomial $T^rw_i$ is of degree $p$. By the definition of
$\prec_{_{T\hbox{-}gr}}$, the definitions of non-central
homogenization and the definition of non-central dehomogenization,
it is not difficult to see that $H$ has the desired property.\par
(iv) By (iii), $F=L+T^r(F_{\sim})^{\sim}$ with $L\in \mathscr{C}$
and $r\in\NZ$. Since by (ii) we have  $F_{\sim}\in \langle~\T
I~\rangle_{\sim}=I$, thus $f=F_{\sim}$ is the desired element. \par
(v) Using basic properties of homogeneous element and graded ideal
in a graded ring, this follows from the foregoing (iii). \QED} \v5

As in the case using central (de)homogenization, before turning to
deal with Gr\"obner bases, we are also concerned about the behavior
of leading monomials under taking the $\NZ$-leading homogeneous
element and non-central (de)homogenization. Below we use $\LH (~)$
to denote taking the $\NZ$-leading homogeneous element (i.e., the
highest-degree homogeneous component) of elements in both $\KS$ and
$\KXT$ with respect to the fixed $\NZ$-gradation.{\parindent=0pt\v5

{\bf 2.2. Lemma}  With the assumptions and notations as fixed above,
the following statements hold.\par (i) If $f\in \KS$, then
$$\LM (f)=\LM (\LH (f))~\hbox{w.r.t.}~\prec_{gr}~\hbox{on}~\B ;$$ If
$F\in\KXT$, then
$$\LM (F)=\LM (\LH
(F))~\hbox{w.r.t.}~\prec_{_{T\hbox{-}gr}}~\hbox{on}~\T{\B}.$$ (ii)
For each nonzero $f\in \KS$, we have
$$\LM (f)=\LM (\T f)~\hbox{w.r.t.}~\prec_{_{T\hbox{-}gr}}~\hbox{on}~\T{\B}.$$
(iii) If $F$ is a homogeneous element in $\KXT$ such that $X_iT{\not
|}~\LM (F)$ with respect to $\prec_{_{T\hbox{-}gr}}$ for all $1\le
i\le n$, then $\LM (F)=T^rw$ for some $r\in\NZ$ and $w\in\B$, such
that $$\LM (F_{\sim})=w=\LM
(F)_{\sim}~\hbox{w.r.t.}~\prec_{gr}~\hbox{on}~\B.$$\par {\bf Proof}
The proof of (i) and (ii) is an easy exercise. To prove (iii), let
$F\in \KXT_p$ be a nonzero homogeneous element of degree $p$. Then
by the assumption $F$ may be written as
$$F=\lambda T^rw+\lambda_1T^{r_1}X_{j_1}W_1+\lambda_2T^{r_2}X_{j_2}W_2\cdots
+\lambda_sT^{r_s}X_{j_s}W_s,$$ where $\lambda,\lambda_i,\in K^*$,
$r,r_i\in\NZ$,  $w\in\B$ and $W_i\in\T{\B}$, such that $\LM
(F)=T^rw$. Since $\B$ consists of $\NZ$-homogeneous elements and the
$\NZ$-gradation of $\KS$ extends to give the $\NZ$-gradation of
$\KXT$, we have $d(T^rw)=d(T^{r_i}X_{j_i}W_i)=p$, $1\le i\le s$.
Also note that $T$ has degree 1. Thus $w=X_{j_i}W_i$ will imply
$r=r_i$ and thereby $T^rw=T^{r_i}X_{j_i}W_i$. So we may assume that
$w\ne X_{j_i}W_i$, $1\le i\le s$. Then it follows from the
definition of $\prec_{_{T\hbox{-}gr}}$ that $r\le r_i$, $1\le i\le
n$. Hence $X_{j_i}W_i\prec_{_{T\hbox{-}gr}}w$, $1\le i\le s$.
Therefore $(X_{j_i}W_i)_{\sim}\prec_{gr}w$, $1\le i\le n$, and
consequently $\LM (F_{\sim})=w=\LM (F)_{\sim}$, as desired.\QED}\v5

The next result strengthens ([Li2], Theorem 8.2), in particular, the
proof of (i) $\Rightarrow$ (ii) given below improves the argument
given in loc. cit. {\parindent=0pt\v5

{\bf 2.3. Theorem} With the notions and notations as fixed above,
let $I=\langle\G\rangle$ be the ideal of $\KS$ generated by a subset
$\G$, and $\langle ~\T I~\rangle$ the non-central homogenization
ideal of $I$ in $\KXT$ with respect to $T$. The following two
statements are equivalent.\par (i) $\G$ is a Gr\"obner basis of $I$
with respect to the admissible system $(\B ,\prec_{gr})$ of
$\KS$;\par (ii) $\T{\G}=\{ \T g~|~g\in \G\}\cup\{X_iT-TX_i~|~1\le
i\le n\}$ is a homogeneous Gr\"obner basis for $\langle~\T
I~\rangle$ with respect to the admissible system
$(\T{\B},\prec_{_{T\hbox{-}gr}})$ of $\KXT$.\vskip 6pt {\bf Proof}
In proving the equivalence below, without specific indication we
shall use (i) and (ii) of Lemma 2.2 wherever it is needed.\par (i)
$\Rightarrow$ (ii) Suppose that $\G$ is a Gr\"obner basis for $I$
with respect to the data $(\B ,\prec_{gr})$. Let $F\in\langle ~\T
I~\rangle$. Then since $\prec_{_{T\hbox{-}gr}}$ is a graded monomial
ordering and hence $\LM (F)=\LM (\LH (F))$, we may assume that $F$
is a nonzero homogeneous element. We want to show that there is some
$D\in\T{\G}$ such that $\LM (D)|\LM (F)$, and hence $\T{\G}$ is a
Gr\"obner basis.}\par Note that $\{ X_iT-TX_i~|~1\le i\le
n\}\subset\T{\G}$ with $\LM (X_iT-TX_i)=X_iT$. If $X_iT|\LM (F)$ for
some $X_iT$, then we are done. Otherwise, $X_iT{\not |}~\LM (F)$ for
all $1\le i\le n$. Thus, by Lemma 2.2(iii), $\LM (F)=T^rw$ for some
$r\in\NZ$ and $w\in\B$ and
$$\LM (F_{\sim})=w=\LM (F)_{\sim}.\leqno{(1)}$$
On the other hand, by Lemma 2.1(iv) we have $F=L+T^{q}\T f$, where
$L$ is an element in the ideal $\mathscr{C}$ generated by $\{
X_iT-TX_i~|~1\le i\le n\}$ in $\KXT$, $q\in\NZ$, and $f\in I$. It
turns out that
$$F_{\sim}=(\T f)_{\sim}=f~\hbox{and hence}~\LM (F_{\sim})=\LM (f).\leqno{(2)}$$
Since $\G$ is a Gr\"obner basis for $I$, there is some $g\in\G$ such
that $\LM (g)|\LM (f)$, i.e., there are $u,v\in\B$ such that
$$\LM (f)=u\LM (g)v=u\LM (\T g)v.\leqno{(3)}$$
Combining (1), (2), and (3) above, we have
$$w=\LM (F_{\sim})=\LM (f)=u\LM (\T g)v.$$
Therefore, $\LM (\T g)|T^rw$, i.e., $\LM (\T g)|\LM (F)$, as
desired.{\parindent=0pt\par (ii) $\Rightarrow$ (i) Suppose that
$\T{\G}$ is a Gr\"obner basis of the graded ideal $\langle~\T
I~\rangle$ in $\KXT$. If $f\in I$, then since $\T f\in \T I$, there
is some $H\in\T{\G}$ such that $\LM (H)|\LM (\T f)$. Note that $\LM
(\T f)=\LM (f)$ and thus $T{\not |}~\LM (\T f)$. Hence $H=\T g$ for
some $g\in\G$, and there are $w,v\in\B$ such that
$$\LM (f)=\LM (\T f)=w\LM (\T g)v=w\LM (g)v.$$
This shows that $\G$ is a Gr\"obner basis for $I$ in $R$. \QED}\v5

We call the Gr\"obner basis $\T{\G}$ obtained in Theorem 2.3 the
{\it non-central homogenization of $\G$} in $\KXT$ with respect to
$T$, or $\T{\G}$ is a {\it non-central homogenized Gr\"obner basis}
with respect to $T$.\v5

By Lemma 2.1 and Theorem 2.3, the following Corollary is
straightforward.{\parindent=0pt\v5

{\bf 2.4. Corollary} Let $I$ be an arbitrary ideal of $\KS$.  With
notation as before, if $\G$ is a Gr\"obner basis of $I$ with respect
to the data $(\B ,\prec_{gr})$, then, with respect to the data
$(\T{\B}, \prec_{_{T\hbox{-}gr}})$ we have
$$\T{\B}-\langle\LM (\T{\G})\rangle =\{ T^rw~|~w\in\B -\langle\LM (\G )\rangle,~r\in\NZ\} ,$$
that is, the set $N(\langle~\T I~\rangle )$ of normal monomials (mod
$\langle~\T I~\rangle )$ in $\T{\B}$ is determined by the set $N(I)$
of normal monomials (mod $I$) in $\B$. Hence, the  algebra
$\KXT/\langle~\T I~\rangle =\KXT /\langle\T{\G}\rangle$ has the
$K$-basis
$$\OV{N(\langle~\T I~\rangle )}=\left\{\left. \OV{T^rw}~\right |~w\in
N(I),~r\in\NZ\right\}.$$}\QED\v5

As with the central (de)homogenization with respect to the commuting
variable $t$ in section 1, theoretically we may also obtain a
Gr\"obner basis for an ideal $I$ of $\KS$ by dehomogenizing a
homogeneous Gr\"obner basis of the ideal $\langle~\T
I~\rangle\subset \KXT$. Below we give a more general approach to
this assertion.{\parindent=0pt \v5

{\bf 2.5. Theorem} Let $J$ be a graded ideal of $\KXT$, and suppose
that $\{ X_iT-TX_i~|~1\le i\le n\}\subset J$. If $\mathscr{G}$ is a
homogeneous Gr\"obner basis of $J$ with respect to the data
$(\T{\B}, \prec_{_{T\hbox{-}gr}})$, then $\mathscr{G}_{\sim}=\{
G_{\sim}~|~G\in\mathscr{G}\}$ is a Gr\"obner basis for the ideal
$J_{\sim}$ in $\KS$ with respect to the data $(\B
,\prec_{gr})$.\vskip 6pt {\bf Proof} If $\mathscr{G}$ is a Gr\"obner
basis of $J$, then $\mathscr{G}$ generates $J$ and hence
$\mathscr{G}_{\sim}=\phi (\mathscr{G} )$ generates $J_{\sim}=\phi
(J)$. We show next that for each nonzero $h\in J_{\sim}$, there is
some $G_{\sim}\in \mathscr{G}_{\sim}$ such that $\LM (G_{\sim})|\LM
(h)$, and hence $\mathscr{G}_{\sim}$ is a Gr\"obner basis for
$J_{\sim}$.}\par Since $\{ X_iT-TX_i~|~1\le i\le n\}\subset J$, by
Lemma 2.1(v) there exists a homogeneous element $H\in J$ and some
$r\in\NZ$ such that $T^r(H_{\sim})^{\sim}=H$ and $H_{\sim}=h$. It
follows that
$$\LM (H)=T^r\LM ((H_{\sim})^{\sim})=T^r\LM (\T h)=T^r\LM (h).\leqno{(1)}$$
On the other hand, there is some $G\in\mathscr{G}$ such that $\LM
(G)|\LM (H)$, i.e., there are $W,V\in\T{\\B}$ such that
$$\LM (H)=W\LM (G)V.\leqno{(2)}$$
But by the above (1) we must have $\LM (G)=T^qw$ for some $q\in\NZ$
and $w\in\B$. Thus, by Lemma 2.2(iii),
$$\LM (G_{\sim})=w=\LM
(G)_{\sim}~\hbox{w.r.t.}~\prec_{gr}~\hbox{on}~\B.\leqno{(3)}$$
Combining (1), (2), and (3) above, we then obtain
$$\begin{array}{rcl} \LM (h)&=&\LM (H)_{\sim}\\
&=&(W\LM (G)V)_{\sim}\\
&=&W_{\sim}\LM (G)_{\sim}V_{\sim}\\
&=&W_{\sim}\LM (G_{\sim})V_{\sim}.\end{array}$$ This shows that $\LM
(G_{\sim})|\LM (h)$, as expected.\QED\v5

We call the Gr\"obner basis $\mathscr{G}_{\sim}$ obtained in Theorem
2.5 the {\it non-central dehomogenization of} $\mathscr{G}$ in $\KS$
with respect to $T$, or $\mathscr{G}_{\sim}$ is a {\it non-central
dehomogenized Gr\"obner basis} with respect to
$T$.{\parindent=0pt\v5

{\bf 2.6. Corollary} Let $I$ be an ideal of $\KS$. If $\mathscr{G}$
is a homogeneous Gr\"obner basis of $\langle~\T I~\rangle$ in $\KXT$
with respect to the data $(\T{\B}, \prec_{_{T\hbox{-}gr}})$, then
$\mathscr{G}_{\sim}=\{ g_{\sim}~|~g\in\mathscr{G}\}$ is a Gr\"obner
basis for $I$ in $\KS$ with respect to the data $(B ,\prec_{gr})$.
Moreover,  if $I$ is generated by the subset $F$ and $\T
F\subset\mathscr{G}$, then $F\subset\mathscr{G}_{\sim}$. \vskip 6pt
{\bf Proof} Put $J=\langle~\T I~\rangle$. Then since $J_{\sim}=I$,
it follows from Theorem 2.5 that if $\mathscr{G}$ is a homogeneous
Gr\"obner basis of $J$ then $\mathscr{G}_{\sim}$ is a Gr\"obner
basis for $I$. The second assertion of the theorem is clear by Lemma
2.1(ii).\QED}\v5

Let $S$ be a nonempty subset of $\KS$ and $I=\langle S\rangle$ the
ideal generated by $S$. Then, with $\T S=\{ \T f~|~f\in S\}\cup\{
X_iT-TX_i~|~1\le i\le n\}$, in general $\langle~\T
S~\rangle\subsetneq\langle~\T I~\rangle$ in $\KXT$ (for instance,
consider $S=\{ Y^3-XY-X-Y,~Y^2-X+3\}$ in the free algebra $K\langle
X,Y\rangle$ and the homogenization in $K\langle X,Y,T\rangle$ with
respect to $T$). Again, as we did in the case dealing with
(de)homogenized Gr\"obner bases with respect to the commuting
variable $t$, we take this place to set up the procedure of getting
a Gr\"obner basis for $I$ and hence a Gr\"obner basis for $\langle
~\T I~\rangle$ by producing a homogeneous Gr\"obner basis of the
graded ideal $\langle~\T S~\rangle$.{\parindent=0pt\v5

{\bf 2.7. Proposition} Let $I=\langle S\rangle$ be the ideal of
$\KS$ as fixed above. Suppose the ground field $K$ is computable.
Then a Gr\"obner basis for $I$ and a homogeneous Gr\"obner basis for
$\langle~\T I~\rangle$ may be obtained by implementing the following
procedure:\par {\bf Step 1.}  Starting with the initial subset $$\T
S=\{ \T f~|~f\in S\}\cup\{ X_iT-TX_i~|~1\le i\le n\} ,$$ compute a
homogeneous Gr\"obner basis $\mathscr{G}$ for the graded ideal
$\langle~\T S~\rangle$ of $\KXT$.\par {\bf Step 2.} Noticing
$\langle~\T S~\rangle_{\sim}=I$, use Theorem 2.5 and dehomogenize
$\mathscr{G}$ with respect to $T$ in order to obtain the Gr\"obner
basis $\mathscr{G}_{\sim}$ for $I$.\par {\bf Step 3.} Use Theorem
2.3 and homogenize $\mathscr{G}_{\sim}$ with respect to $T$ in order
to obtain the homogeneous Gr\"obner basis
$(\mathscr{G}_{\sim})^{\sim}$ for the graded ideal $\langle~\T
I~\rangle$.\par\QED}\v5

\v5 \centerline{References}
\parindent=1truecm

\item{[Bu]} B. Buchberger, Gr\"obner bases: An algorithmic method in
polynomial ideal theory. In: {\it Multidimensional Systems Theory}
(Bose, N.K., ed.), Reidel Dordrecht, 1985, 184--232.

~\re{[BW]} T.~Becker and V.~Weispfenning, {\it Gr\"obner Bases},
Springer-Verlag, 1993.

\item{[Gr]} E. L. Green, Noncommutative Gr¡§obner bases and projective
resolutions, in: {\it Proceedings of the Euroconference
Computational Methods for Representations of Groups and Algebras},
Essen, 1997, (Michler, Schneider, eds), Progress in Mathematics,
Vol. 173, Basel, Birkha¡§user Verlag, 1999, 29--60.

\item{[HT]} D. Hartley and P. Tuckey, Gr\"obner Bases in Clifford and
Grassmann Algebras, {\it J. Symb. Comput.}, 20(1995), 197--205.

\item{[KRW]} A.~Kandri-Rody and V.~Weispfenning, Non-commutative
Gr\"obner bases in algebras of solvable type, {\it J. Symbolic
Comput.}, 9(1990), 1--26.

\item{[Li1]} H. Li, {\it Noncommutative Gr\"obner Bases and
Filtered-Graded Transfer}, LNM, 1795, Springer-Verlag, 2002.

\item{[Li2]}  H. Li, $\Gamma$-leading homogeneous algebras and Gr\"obner
bases, {\it Advanced Lectures in Mathematics}, Vol.8, International
Press and Higher Education Press, Boston-Beijing, 2009, 155--200.

\item{[Li3]} H. Li, On the calculation of gl.dim$G^{\NZ}(A)$ and gl.dim$\widetilde{A}$
 by using Groebner bases, {\it Algebra Colloquium}, 16(2)(2009), 181--194.

\item{[LWZ]} H. Li, Y. Wu and J. Zhang, Two applications of
noncommutative Gr\"obner bases, {\it Ann. Univ. Ferrara - Sez. VII -
Sc. Mat.}, XLV(1999), 1--24.

\item{[Mor]} T. Mora, An introduction to commutative and noncommutative
Gr\"obner Bases, {\it Theoretic Computer Science}, 134(1994),
131--173.

\item{[SL]} C. Su and H. Li, Some results on dh-closed homogeneous
Gr\"obner bases and dh-closed graded ideals, arXiv:math.RA/0907.0526

\end{document}